\documentclass{article}
\usepackage{amsmath}
\usepackage{amssymb}
\usepackage{graphicx}
\usepackage{breqn}
\usepackage{epstopdf}
\title{A note on maximum likelihood estimation for cubic and quartic canonical toric Del Pezzo surfaces}
\author{Dimitra Kosta}
\date{}

\newtheorem{theorem}{Theorem}[section]

\newtheorem{proposition}[theorem]{Proposition}

\newtheorem{definition}[theorem]{Definition}

\begin{document}
 \maketitle

\begin{abstract}
This article focuses on the study of algebraic statistical models corresponding to toric Del Pezzo surfaces with Du Val singularities. A closed-form for the maximum likelihood estimate of algebraic statistical models which correspond to cubic and quartic toric Del Pezzo surfaces with Du Val singular points is given. Also, we calculate the ML degrees of some toric Del Pezzo surfaces of degree less than or equal to six and show that they equal the degree of the surface in all the cases but one, namely the quintic Del Pezzo with two points of type $\mathbb{A}_1$.
\end{abstract}

\section{Introduction}
Maximum likelihood estimation (MLE) is a standard approach to parameter estimation and inference, and a fundamental computational task in statistics. It consists of the following problem: given the observed data and a model of interest, find the probability distribution that is most likely to have produced the data. In the past decade, algebraic techniques for the computation of maximum likelihood estimates have been developed for algebraic statistical models for discrete data (see \cite{Catanese_etal},  \cite{Hauenstein_etal}, \cite{Hosten_etal}).
\let\thefootnote\relax\footnote{The author would like to thank Ivan Cheltsov, Alexander Davie, Milena Hering, Kaie Kubjas and Bernd Sturmfels for valuable comments and suggestions. This article was completed during a Daphne Jackson Trust Fellowship, funded jointly by EPSRC and the University of Edinburgh.}

These efforts to solve the likelihood equations for algebraic statistical models have focused on specifying the ML-degree of the corresponding models (see \cite{Catanese_etal},  \cite{Hauenstein_etal}, \cite{Hosten_etal} and \cite{Lin_etal}). In \cite{Hosten_etal} algebraic algorithms are presented for computing all critical points of the likelihood function, with the aim of identifying the local maxima in the probability simplex. Using these algorithms the authors derived an explicit formula for the ML-degree of a projective variety which is a generic complete intersection and this formula serves as an upper bound for the ML-degree of special complete intersections. Moreover, a geometric characterisation of the ML-degree of a smooth variety in the case when the divisor corresponding to the rational function is a normal crossings divisor was given in \cite{Catanese_etal}. In the same paper an explicit combinatorial formula for the ML-degree of a toric variety is derived by relaxing the restrictive smoothness assumption and allowing some mild singularities. Also, in \cite{Huh} a geometric description of all varieties with ML-degree 1 is given via Horn Uniformization, however given a variety it is hard to tell if that variety has ML-degree 1. For an introduction to the geometry behind the MLE for algebraic statistical models for discrete data the interested reader is refered to \cite{Huh_Sturmfels}, which includes most of the current results on the MLE problem from the perspective of algebraic geometry as well as statistical motivation.

This article is concerned with the problem of MLE  for algebraic statistical models with singularities, in particular those which correspond to toric Del Pezzo surfaces with Du Val singular points (also called canonical toric Del Pezzo surfaces). The focus is not towards specifying the ML degree, but towards providing a closed form for the maximum likelihood estimate. Throughout this paper, we work over the field $\mathbb{C}$ of complex numbers. In particular we prove the following result.

\begin{theorem}
\label{MainResult}
Let $X_d$ be a toric Del Pezzo surface of degree $d$ with Du Val singularities. Then there is a closed form for the unique maximum likelihood
estimate for the corresponding algebraic statistical model only in the case of a cubic $X_3$ and a quartic $X_4$ toric Del Pezzo surface.
\end{theorem}

The importance of algebraic statistical models corresponding to toric varieties is  due to their correspondence to
log-linear models in Statistics (see \cite{Drton_etal}). Since the seminal papers by L.A. Goodman in the 70s, log-linear
models have been widely used in statistics and areas like natural language processing when analysing crossclassified
data in multidimensional contingency tables \cite{Bishop_etal}.

Another reason for studying the MLE for such algebraic statistical models, is that they correspond to singular varieties.
Singularities play an important role in statistical inference as the commonly assumed smoothness of algebraic statistical models is very restrictive and is almost never satisfied for models of statistical relevance (see  \cite[p. 100]{Drton_etal}, \cite{Sturmfels2009}). For instance, algebraic
statistical models of binary symmetric phylogenetic three-valent trees are proven to be Fano varieties with Gorenstein terminal singularities in \cite{Wisniewski_Buczynska_2007}. Statistical applications point towards a better understanding of the MLE problem for higher dimensional singular toric Fano models and  studying the MLE problem for algebraic statistical models corresponding to toric Del Pezzo surfaces with Du Val singular points is the natural first step.

\bigskip
\section{Preliminaries}
The relevant definitions to this problem are given in this section. In section \ref{MLE} we  introduce the reader to the MLE problem from the perspective of algebraic geometry. In section \ref{toric}, we define toric models following the notation used in Chapter 1.2 of \cite{Pachter_Sturmfels}. Section \ref{delPezzo} serves as a brief introduction to canonical Del Pezzo surfaces with a focus on the properties of the toric ones that will be used in the proof of Theorem~\ref{MainResult}.




\subsection{Maximum Likelihood Estimation}
\label{MLE}

\bigskip

Consider the complex projective space $\mathbb{P}^n$ with coordinates $(p_0, p_1,...,p_n)$. In our setting we will consider $X$ to be a discrete random variable taking values on the state space $[n]$. Then the coordinate $p_i$ represents the probability of the $i$-event $$p_i = P (X = i) \text{,} $$ where $i=0, \cdots, n$, therefore $p_0+p_1+...+p_n =1$. The set of points in $\mathbb{P}^n$ with positive real coefficients is identified with the probability simplex
$$
\Delta_n =\{  (p_0, p_1, ... , p_n) \in \mathbb{R}^{n+1} : p_0, p_1, ... , p_n \geq 0 \text{ and } p_0+p_1+...+p_n =1 \} \text{ . }
$$
An algebraic statistical model is a Zariski closed subset $\mathcal{M}$ of complex projective space  $\mathbb{P}^n$, with the model itself being the intersection of $\mathcal{M}$  with the probability simplex $\Delta_n$. The data is given by a non-negative integer vector $(u_0, u_1,...,u_n) \in \mathbb{N}^{n+1}$, where $u_i$ is the number of times the $i$-event is observed.

The maximum likelihood estimation problem consists of finding a model point $p \in \mathcal{M}$ which maximises the likelihood of observing the data. This amounts to maximising the corresponding likelihood function
$$
L(p_0, p_1,...,p_n) = \frac{p_0^{u_0} \cdot p_1^{u_1} \cdots p_n^{u_n}}{(p_0+p_1+...+p_n)^{(u_0+u_1+...+u_n)}}
$$
over the model $\mathcal{M}$, where here we ignore a multinomial coefficient. Statistical computations are usually implemented in the affine $n$-plane   $p_0+p_1+...+p_n =1$. However, including the denominator makes the likelihood function a well-defined rational function on the projective space  $\mathbb{P}^n$, enabling one to use projective algebraic geometry to study its restriction to the variety $\mathcal{M}$.

The likelihood function might not be convex, so it can have many local maxima and the problem of finding and certifying a global maximum is difficult. Therefore, in most recent works the problem of finding all critical points of the likelihood function is considered, with the aim of identifying all local maxima (see \cite{Catanese_etal}, \cite{Hauenstein_etal} and \cite{Hosten_etal}). This corresponds to solving a system of polynomial equations and these equations, defining the critical points of the likelihood function $L$, are called likelihood equations. The number of complex solutions to the likelihood equations equals the number of complex critical points of the restriction of the likelihood function $L$ to the model $\mathcal{M}$, which is called the maximum likelihood (ML) degree of the variety $\mathcal{M}$.

\subsection{Toric models}
\label{toric}

\bigskip

In this article we are studying the maximum likelihood estimation problem for toric models which are models with a well behaved likelihood function.  Toric models are known as log-linear models in statistics, because the logarithms of the probabilities are linear functions in the
logarithms of the parameters $\theta_i$.  They have the property that maximum likelihood estimation  is a convex optimization problem. Assuming that the parameter domain $\Theta$ is bounded, it follows that the likelihood function has exactly one local maximum.

Let $A = ( a_{ij} )$ be a non-negative
integer $d \times m$ matrix with the property that all column sums are equal:
$$
\sum_{i=1}^{d} a_{i1} = \sum_{i=1}^{d} a_{i2} = ... = \sum_{i=1}^{d} a_{im} \text{ . }
$$
The $j$-th column vector $a_j$ of the matrix $A$ represents the monomial
$$
\theta^{a_j} := \theta_1^{a_{1j}} \cdot \theta_2^{a_{2j}} \cdot \cdot \cdot \theta_d^{a_{dj}}  \text{  for all  }  j=1,...m
$$
and the assumption that the column sums of the matrix $A$ are all equal means these monomials all have the same degree.

\begin{definition}
The toric model of $A$ is the image of the orthant $\Theta = \mathbb{R}^d_{>0}$ under the map
$$
f : \mathbb{R}^d \to \mathbb{R}^m, \theta \mapsto \frac{1}{\sum_{j=1}^{m} \theta^{a_j}} \cdot (\theta^{a_1}, \theta^{a_2}, ..., \theta^{a_m}) \text{ .}
$$

\end{definition}

Maximum likelihood estimation for the toric model means solving the optimization problem of maximising the function
$$
p_1^{u_1} \cdot p_1^{u_2} \cdot ... \cdot p_m^{u_m}
$$
subject to the constrains $f(\mathbb{R}^d_{>0})$. This is equivalent to maximising function
$$
\theta^{A \cdot u}  \text{  subject to  } \theta \in \mathbb{R}^d_{>0} \text{  and  } \sum_{j=1}^{m} \theta^{a_j} =1 ,
$$
where
$$
\theta^{Au} = \prod_{i=1}^{d} \theta_i^{a_{i1}u_1+a_{i2}u_2+...+ a_{im}u_m}  \text{  and  } \theta^{a_{j}} = \prod_{i=1}^{d} \theta_i^{a_{ij}} \text{. }
$$
Let $b:= Au$ denote the sufficient statistic, then the optimisation problem above is equivalent to maximising
$$
 \theta^b \text{ subject to } \theta \in \mathbb{R}^d_{>0}  \text{  and  }  \sum_{j=1}^{m} \theta^{a_j} =1 \text{ .}
$$ Then by the next Proposition \cite[Proposition 1.9]{Pachter_Sturmfels} we get an equation satisfied by the local maximum.
\begin{proposition}

Fix a toric model $A$ and data $u \in \mathbb{N}^m$ with sample size
$N = u_1 + \cdots+ u_m$ and sufficient statistic $b = Au$. Let $\hat{p} = f (\hat{\theta})$ be any local
maximum for the equivalent optimization problems above. Then
$$
A \cdot \hat{p} =  \frac{1}{N} \cdot b \text{.}
$$
\end{proposition}

Given a matrix $A \in \mathbb{N}^{d \times m}$ and any vector $b \in \mathbb{R}^d$, we consider the set
$$
P_A(b) = \{ p \in \mathbb{R}^m : A \cdot p = \frac{1}{N} \cdot b \text{  and  } p_j>0 \text{ for all } j \} .
$$
This is a relatively open polytope and Birch's theorem \cite[Theorem 1.10]{Pachter_Sturmfels} below asserts that it is either empty or meets the toric model in precisely one point.

\begin{theorem}[Birch's Theorem]
Fix a toric model $A$ and let $u \in \mathbb{N}^m_{>0}$ be a strictly positive data vector with sufficient statistic $b = A u$. The intersection of the polytope $P_A(b)$ with the toric model $f(\mathbb{R}^d_{>0})$ consists of precisely one point. That point is the maximum likelihood estimate $\hat{p}$ for the data $u$.
\end{theorem}

\subsection{Del Pezzo surfaces}
\label{delPezzo}

This section aims to provide a brief introduction to Del Pezzo surfaces with the reader being referred to \cite{Dolgachev} for a more detailed study of this classical subject of algebraic geometry. In particular, we are focusing on the characterisation of toric Del Pezzo surfaces with Du Val singularities in terms of their polytopes which will be used in the proof of Theorem \ref{MainResult}.

Del Pezzo surfaces were named after the Italian mathematician Pasquale Del Pezzo who encountered this class of surfaces when studying surfaces of degree $d$ embedded in $\mathbb{P}^d$. They are two dimensional examples of Fano varieties, a class of varieties with ample anticanonical divisor class which have been extensively studied in birational geometry in the context of the minimal model program (see \cite{Kollar}).

\begin{definition}
A complex projective algebraic surface $X$ with ample anticanonical divisor class $-K_X$ is called a Del Pezzo surface.
\end{definition}

The self-intersection number $K^2_X$ of the canonical class of the Del Pezzo surface is called the degree $d$ of the Del Pezzo surface $X$.
Often Del Pezzo surfaces are assumed to be smooth, but here our Del Pezzo surfaces are allowed to have mild singularities and, in particular, at most Du Val singularities.

\begin{definition}
A point $P$ of a normal surface $X$ is called a Du Val singularity if there exists a minimal resolution $\pi :  \tilde{X} \to X$ such that $K_{\tilde{X}} \cdot E_i = 0$ for every exceptional curve $E_i \subset  \tilde{X}$.
\end{definition}

Every Del Pezzo surface is either a product of two projective lines $\mathbb{P}^1 \times \mathbb{P}^1$  (with d=8), or the blow-up of a projective plane in $9-d$ points with no three collinear, no six on a conic, and no eight of them on a cubic having a node at one of them. Conversely any blowup of the plane in points satisfying these conditions is a Del Pezzo surface. In dimension two, the class of Du Val singularities is the same as the class of canonical singularities firstly introduced by M. Reid (see \cite{Reid}), although this is not the case in higher dimensions. This is why we refer to Del Pezzo surfaces with Du Val singularities also as canonical Del Pezzo surfaces.

Toric Del Pezzo surfaces correspond to certain convex lattice polytopes whose boundary lattice points are dictated by the singularities involved. In particular, the corresponing polytope of a toric Del Pezzo surface with Du Val singularities is a reflexive polytope. According to the classification results of reflexive polytopes (see \cite{Kreuzer_Skarke}, \cite{Sato}), we have the following.

\begin{proposition} There are exactly 16 isomorphism classes of two-dimensional reflexive polytopes given in the list bellow. The number in the labels is the number of lattice points on the boundary.
\end{proposition}

\begin{table}[h!]
\centering
 \includegraphics[scale=0.80]{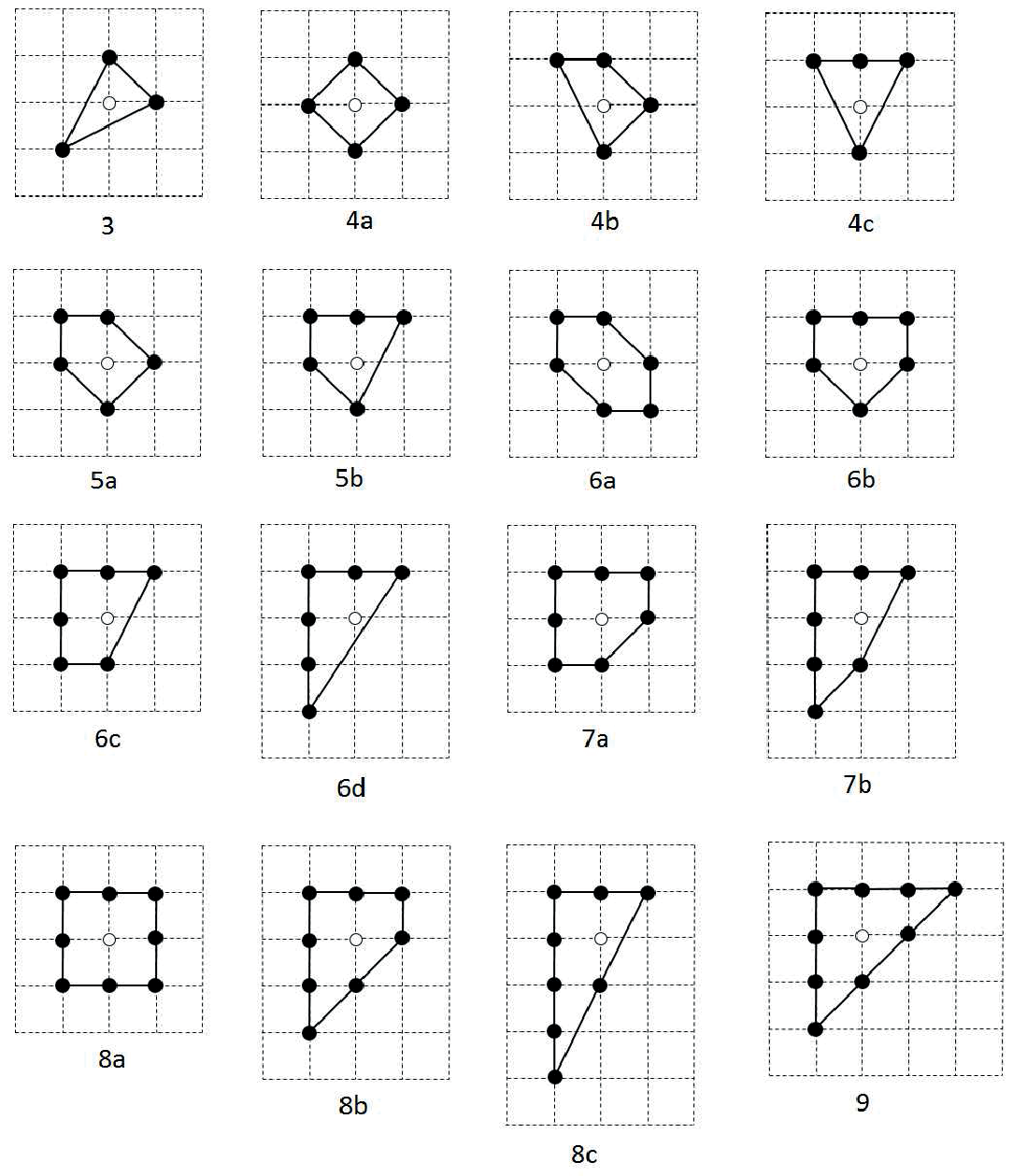}
 \caption{Two-dimensional reflexive polytopes}
\label{table:2}
\end{table}

\section{Main Calculation}
In this section we determine a closed form for the maximum likelihood  estimates of all algebraic statistical models corresponding to cubic (d=3) and quartic (d=4) toric Del Pezzo surfaces  with Du Val singularities. In each case we will explain how one can detect the type and number of singular points from the corresponding polytope.

\subsection{Cubic Del Pezzo with three singular points of type $\mathbb{A}_2$}
Consider the case of a reflexive polytope in $\mathbb{Z}^2$ with three lattice points $(1,0)$, $(0,1)$ and $(-1,-1)$ on the boundary, as in the graph below.
\begin{center}
 \includegraphics[scale=0.35]{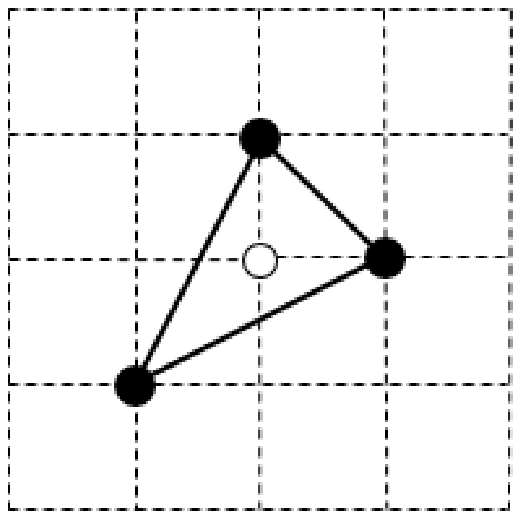}
\end{center}
By computing the normal fan we see its minimal generators span a sublattice of index $3$ in $\mathbb{Z}^2$. This corresponds to a cubic surface with three Du Val singular points of type $\mathbb{A}_2$. Another way to see this is by using the characters coming from the lattice points of the polygon. The corresponding toric variety is isomorphic to the cubic surface with three double points of type $\mathbb{A}_2$ and equation
$S_3:  p_4^3 - p_1 p_2 p_3 = 0$ in $\mathbb{P}^3$.

This polytope generated by the lattice points
 $(1,0),(0,1),(-1,-1)$ in $\mathbb{Z}^2$ gives rise to the same projective toric variety as the polytope which is generated by the  lattice points $ (2,1,0), (1,2,0), (0,0,3), (1,1,1) $ in $\mathbb{Z}^3$.

We are interested in the algebraic statistical model given by the matrix
\[
A=
  \begin{bmatrix}
    2 & 1 & 0 & 1\\
    1 & 2 & 0 & 1\\
    0 & 0 & 3 & 1
  \end{bmatrix} \text{ . }
\]
This non-negative integer matrix $A$ corresponds to the map
 $$f: \mathbb{C}^3 \to \mathbb{C}^3  \text{  ,  }   (\theta_1, \theta_2, \theta_3) \mapsto \frac{(\theta_1^2 \theta_2, \theta_1 \theta_2^2, \theta_3^3, \theta_1 \theta_2 \theta_3 )}{(\theta_1^2 \theta_2 + \theta_1 \theta_2^2 + \theta_3^3 + \theta_1 \theta_2 \theta_3)} \text{ . }$$	

According then to Birch's theorem there is a unique maximum likelihood estimate $\hat{p}$ for the data $u = (u_1, u_2, u_3, u_4)$ with $N =u_1 + u_2 +u_3 + u_4$. This unique MLE satisfies the equation
$$
A \cdot \hat{p} =  \frac{1}{N} \cdot A \cdot u \text{ ,}
$$
where $N =u_1 + u_2 + u_3 + u_4$ and $\hat{p}$ is the probability distribution corresponding to the parameter values $(\hat{\theta}_1, \hat{\theta}_2, \hat{\theta}_3) $. This gives us the equations
\begin{eqnarray*}
\label{eq:likeq_1}
3 \hat{p}_1 + \hat{p}_4 & = & \frac{1}{N} (3 u_1 + u_4) \\
3\hat{p}_2 +  \hat{p}_4 & = & \frac{1}{N} (3 u_2 +  u_4) \\
3 \hat{p}_3 + \hat{p}_4 & = & \frac{1}{N} (3u_3 + u_4) \text{ . }
\end{eqnarray*}
 We then compute $(\hat{\theta}_1, \hat{\theta}_2, \hat{\theta}_3) = (  \sqrt[3]{ \frac{p_1^2}{p_2}} , \sqrt[3]{ \frac{p_2^2}{p_1}} , \sqrt[3]{p_3} )$,
 where each $p_k \text{, } k=1,2,3$ is given by a cubic equation

 \begin{dmath*}
 \hat{p}_k^3 - \frac{N-28(3u_k+u_4)}{28N}\hat{p}_k^2 + \frac{(u_i-u_k)(u_j-u_k)-9(3u_k+u_4)^2}{28N^2}\hat{p}_k  - \frac{(3u_k+u_4)^3}{28N^3} = 0
\end{dmath*}
for $i,j,k = 1,2,3$ and $i \neq j \neq k$.

\subsection{Quartic Del Pezzo with four$\mathbb{A}_1$ type singular points.}
Consider the reflexive polytope with four lattice points on the boundary, as in the graph below.

\begin{center}
 \includegraphics[scale=0.35]{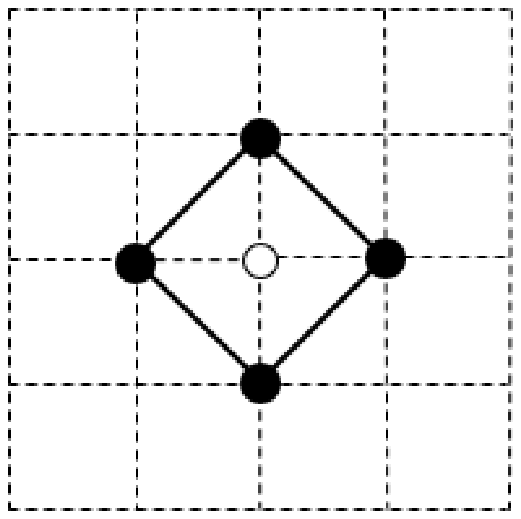}
\end{center}

If we compute the normal fan, we see that it corresponds to a quartic surface with four Du Val singular points of type $\mathbb{A}_1$. Another way to see this, is by using the characters coming from the lattice points of the polygon. Then the corresponding projective toric variety is isomorphic to a quartic surface. This is the complete intersection of the two quadric equations $S_4:  p_1 p_4-p_5^2 = p_2 p_3-p_5^2 = 0$ in $\mathbb{P}^4$, which has four double points of type $\mathbb{A}_1$. The polytope generated by the lattice points $(1,0), (0,1), (0,-1),(-1,0)\in \mathbb{Z}^2$ gives rise to the same projective toric variety as the polytope generated by the lattice points $(2,1,0),(1,2,0),(1,0,2),(0,1,2), (1,1,1)$ in $\mathbb{Z}^3$.

The non-negative integer matrix
\[
A=
  \begin{bmatrix}
    2 & 1 & 1 & 0 & 1 \\
    1 & 2 & 0 & 1 & 1\\
    0 & 0 & 2 & 2 & 1
  \end{bmatrix}
\]
corresponds to the map
$$f: \mathbb{C}^3 \to \mathbb{C}^4  \text{  ,  }   (\theta_1, \theta_2, \theta_3) \mapsto \frac{(\theta_1^2 \theta_2, \theta_1 \theta_2^2, \theta_1 \theta_3^2, \theta_2 \theta_3^2 , \theta_1 \theta_2 \theta_3)}{(\theta_1^2 \theta_2 + \theta_1 \theta_2^2 + \theta_1 \theta_3^2+ \theta_2 \theta_3^2 ) } \text{ . }$$	
According to Birch's theorem there is a unique maximum likelihood estimate $\hat{\theta}$ for the data $u = (u_1, u_2, u_3, u_4, u_5)$ which satisfies the equation
$$
A \cdot \hat{p} =  \frac{1}{N} \cdot A \cdot u \text{ , }
$$
where $N =u_1 + u_2 + u_3 + u_4 +u_5$ and $\hat{p}$ is the probability distribution corresponding to the parameter values $(\hat{\theta}_1, \hat{\theta}_2, \hat{\theta}_3) $.
This gives us the equations
\begin{eqnarray*}
\label{eq:likeq_2}
\hat{p}_1 - \hat{p}_4 & = & \frac{1}{N} (u_1 -u_4) \\
\hat{p}_2 - \hat{p}_3 & = & \frac{1}{N} (u_2 - u_3) \\
2 \hat{p}_3 + 2\hat{p}_4 - \hat{p}_5 & = & \frac{1}{N} (2u_3 + 2u_4 - u_5) \text{. }
\end{eqnarray*}
Then  the unique maximum likelihood estimate $\hat{\theta}$ for the data $u$  is
$(\hat{\theta}_1, \hat{\theta}_2, \hat{\theta}_3) = (\sqrt[3]{ \frac{\hat{p}_1^2}{\hat{p}_2}} , \sqrt[3]{\frac{\hat{p}_2^2}{\hat{p}_1}} , \sqrt[3]{\frac{\hat{p}_5^3} {\hat{p}_1 \hat{p}_2}})$,
 where $\hat{p}_1, \hat{p}_2, \hat{p}_5$ are given by the quartic equations below
 \begin{dmath*}
 \Big[ \hat{p}_1^2 - \frac{4N+(u_1 -u_4)}{N}\hat{p}_1 + \frac{(2u_1+2u_2+u_5)(2u_1+2u_3+u_5)}{N^2} \Big]^2 - \Big[\hat{p}_1^2-\frac{(u_1-u_4)\hat{p}_1}{N} \Big] \Big[4\hat{p}_1 - \frac{2N +2(u_1-u_4)}{N} \Big]^2 = 0 \end{dmath*}
 and
\begin{dmath*}
 \Big[ \hat{p}_2^2 - \frac{4N+(u_2 -u_3)}{N}\hat{p}_2 + \frac{(2u_1+2u_2+u_5)(2u_2+2u_4+u_5)}{N^2} \Big]^2 - \Big[\hat{p}_2^2-\frac{(u_2-u_3)\hat{p}_2}{N} \Big] \Big[4\hat{p}_2 - \frac{2N +2(u_2-u_3)}{N} \Big]^2 = 0 \end{dmath*}
 and
\begin{dmath*} \Big[\hat{p}_5^2 - 2\hat{p}_5+\frac{(2u_1+2u_2+u_5)(2u_3+2u_4+u_5)}{N^2}+\frac{(2u_1-2u_4)(N+u_2-u_3)}{N^2} \Big]^2 - (2-2\hat{p}_5)^2 \big[ \frac{(u_1-u_4)}{N^2} +4\hat{p}_5^2 \big] = 0 \end{dmath*}

\subsection{Quartic Del Pezzo with one $\mathbb{A}_2$ and two $\mathbb{A}_1$ type singular points}
Consider the case of a reflexive polytope with four lattice points on the boundary. By using the characters coming from the lattice points of the polygon, the corresponding toric variety is isomorphic to a quartic surface. This is the complete intersection of two quadrics $S_4^{''}: p_1 p_3 - p_2 p_5= p_2 p_4 -p_5^2 = 0$ in $\mathbb{P}^4$, which has one $\mathbb{A}_2$ and two $\mathbb{A}_1$ type singular points.

The polytope generated by the lattice points $(1,0), (0,1), (-1,1), (0,-1)$ in $\mathbb{Z}^2$ gives rise to the same projective toric variety as the polytope generated by the lattice points $ (2,1,0), (1,2,0), (0,2,1), (1,0,2),  (1,1,1)$ in $\mathbb{Z}^3$.

\begin{center}
 \includegraphics[scale=0.35]{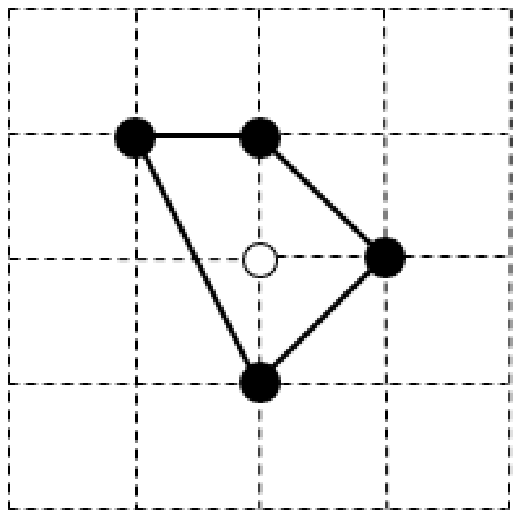}
\end{center}

The non-negative integer matrix
\[
A=
  \begin{bmatrix}
   2 & 1 & 0 & 1 & 1\\
   1 & 2 & 2 & 0 & 1\\
   0 & 0 & 1 & 2 & 1
  \end{bmatrix}
\]
corresponds to the map
$$f: \mathbb{C}^3 \to \mathbb{C}^4  \text{  ,  }   (\theta_1, \theta_2, \theta_3) \mapsto \frac{( \theta_1^2 \theta_2, \theta_1 \theta_2^2, \theta_2^2 \theta_3, \theta_1 \theta_3^2, \theta_1 \theta_2 \theta_3 )}{( \theta_1 \theta_3^2 + \theta_2^2 \theta_3 + \theta_1 \theta_2^2 + \theta_1^2 \theta_2 + \theta_1 \theta_2 \theta_3) } \text{ . }$$	

According to Birch's theorem there is a unique maximum likelihood estimate $\hat{\theta}$ for the data $u = (u_1, u_2, u_3, u_4, u_5)$ which satisfies the equation
$$
A \cdot \hat{p} =  \frac{1}{N} \cdot A \cdot u \text{ , }
$$
where $N =u_1 + u_2 + u_3 + u_4+u_5$ and $\hat{p}$ is the probability distribution corresponding to the parameter values $(\hat{\theta}_1, \hat{\theta}_2, \hat{\theta}_3) $.

This gives us the equations
\begin{eqnarray*}
\label{eq:likeq_3}
\hat{p}_1 + 2\hat{p}_4 + \hat{p}_5 & = & \frac{1}{N} (u_1 + 2u_4 + u_5) \\
\hat{p}_2 -  3\hat{p}_4 - \hat{p}_5 & = & \frac{1}{N} (u_2 -  3u_4 - u_5) \\
\hat{p}_3 + 2\hat{p}_4 + \hat{p}_5 & = & \frac{1}{N} (u_3 + 2u_4 + u_5)
 \text{ . }
\end{eqnarray*}
Then  the unique maximum likelihood estimate $\hat{\theta}$ for the data $u$  is
$(\hat{\theta}_1, \hat{\theta}_2, \hat{\theta}_3) = (\sqrt[3]{ \frac{\hat{p}_1^2}{\hat{p}_2}} , \sqrt[3]{\frac{\hat{p}_2^2}{\hat{p}_1}} , \sqrt[3]{\frac{\hat{p}_5^3} {\hat{p}_1 \hat{p}_2}})$,
 where $\hat{p}_1$ is given by the quartic equation below
  \begin{dmath*}
 \Big[ 9\hat{p}_2^2 + \frac{4(3u_1+2u_2+u_5)-6(u_1 -u_3)}{N}\hat{p}_2 + \frac{(u_1-u_3)^2}{N^2} \Big] \Big[
 16 \hat{p}_2 + \frac{6(3u_1+2u_2+u_5)-9(u_1 -u_3)}{N} \Big]^2 - \Big[33\hat{p}_2^2+\frac{23(3u_1+2u_2+u_5)-43(u_1 -u_3) +(u_1+u_2-u_4)}{N}\hat{p}_2 + \frac{(3u_3+2u_2+u_5)^2}{N^2} \Big]^2 = 0 \end{dmath*}
$\hat{p}_2$ is given by the quartic equation
   \begin{dmath*}
 \Big[ \hat{p}_1^2 + \frac{8(u_1 -u_3)-6(3u_1+2u_2+u_5)}{N}\hat{p}_1 + \frac{(3u_1+2u_2+u_5)^2}{N^2} \Big] \Big[
 19 \hat{p}_1 + \frac{2(u_1+u_2-u_4)-7(3u_1+2u_2+u_5)}{N} \Big]^2 - \Big[3\hat{p}_1^2+\frac{12(u_1 -u_3)-8(3u_1+2u_2+u_5)-6(u_1+u_2-u_4)}{N}\hat{p}_1 + \frac{(3u_3+2u_2+u_5)^2 +2(3u_3+2u_2+u_5)(u_1+u_2-u_4)}{N^2} \Big]^2 = 0 \end{dmath*}
 and the probability distribution $ \hat{p}_5$ is given by the quartic equation
$$
a_5 \hat{p}_5^4 + b_5 \hat{p}_5^3 + c_5 \hat{p}_5^2 + d_5 \hat{p}_5 + e_5 =0
$$
where
\begin{eqnarray*}
a_5 & = & 51  \\
b_5 & = & -35(u+v+w)-4v  \\
c_4 & = & (u+v+w)(9u+36v+w)+75uw-8v^2 \\
d_4 & = &  -3(u+v+w)\big[ 12uw+(6v+6w+4v)v \big] -3uvw\\
e_4 & = & 3uw\big[ 9uw+(6u+6w+4v)v \big]  \text{ , }
\end{eqnarray*}
with $u=\frac{u_1+2u_4+u_5}{N}$, $v=\frac{u_3+2u_4+u_5}{N}$ and $w=\frac{u_2-3u_4-u_5}{N}$.

\subsection{Quartic Del Pezzo surface with one $\mathbb{A}_3$ and two $\mathbb{A}_1$ type singular points.}
Consider the case of a reflexive polytope with four lattice points on the boundary, as in the graph below.

\begin{center}
 \includegraphics[scale=0.35]{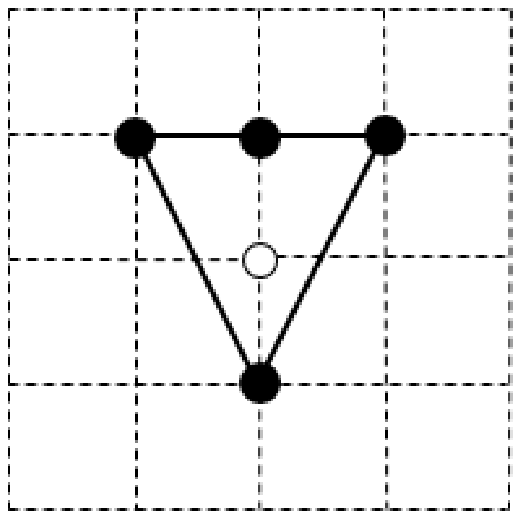}
\end{center}
 We use the characters coming from the lattice points of the polygon. Then the corresponding toric variety is isomorphic to a quartic surface, which is the complete intersection of two quadrics $S_4^{'''}: p_2p_4-p_3^2 = p_1p_3 - p_5^2 = 0$ in $\mathbb{P}^4$, which has one $\mathbb{A}_3$ and two $\mathbb{A}_1$ type singular points.  We can see that the polytope generated by the lattice points $(1,1), (0,1), (-1,1)$ and $(0,-1)$ in $\mathbb{Z}^2$ gives rise to the same projective toric variety as the polytope generated by the lattice points $(1,0,3),(2,2,0),(1,2,1),(0,2,2), (1,1,2)$ in $\mathbb{Z}^3$.

The non-negative integer matrix
\[
A=
  \begin{bmatrix}
    1 & 2 & 1 & 0 & 1 \\
    0 & 2 & 2 & 2 & 1 \\
    3 & 0 & 1 & 2 & 2
  \end{bmatrix}
\]
gives the map
$$f: \mathbb{C}^3 \to \mathbb{C}^4  \text{  ,  }   (\theta_1, \theta_2, \theta_3) \mapsto \frac{(\theta_1 \theta_3^3, \theta_1^2 \theta_2^2, \theta_1 \theta_2^2 \theta_3, \theta_2^2 \theta_3^2, \theta_ 1 \theta_2 \theta_3^2  )}{( \theta_1 \theta_3^3 + \theta_1^2 \theta_2^2+ \theta_1 \theta_2^2 \theta_3 + \theta_2^2 \theta_3^2 ) } \text{ . }$$	

According to Birch's theorem there is a unique maximum likelihood estimate $\hat{\theta}$ for the data $u = (u_1, u_2, u_3, u_4, u_5)$ which satisfies the equation
$$
A \cdot \hat{p} =  \frac{1}{N} \cdot A \cdot u \text{ , }
$$
where $N =u_1 + u_2 + u_3 + u_4 + u_5$ and $\hat{p}$ is the probability distribution corresponding to the parameter values $(\hat{\theta}_1, \hat{\theta}_2, \hat{\theta}_3) $.

This gives us the equations
\begin{eqnarray*}
\label{eq:likeq_4}
\hat{p}_1 - \hat{p}_4 & = & \frac{1}{N} (u_1 - u_4)  \\
2 \hat{p}_2 + \hat{p}_4 + \hat{p}_5 & = & \frac{1}{N} (2u_2 + u_4 + u_5) \\
2 \hat{p}_3 + \hat{p}_5 & = & \frac{1}{N} (2u_3 + u_5) \text{ . }
\end{eqnarray*}
We can compute the unique maximum likelihood estimate $\hat{\theta}$ for the data $u$, which is
\begin{eqnarray*}
\hat{\theta}_1 & = & \frac{\hat{p}_2}{\hat{p}_3} \sqrt[4]{\frac{\hat{p}_5^2}{\hat{p}_2}} \\
\hat{\theta}_2 & = & \hat{p}_3 \sqrt[4]{\frac{\hat{p}_2^3}{\hat{p}_5^6}} \\
\hat{\theta}_3 & = & \sqrt[4]{\frac{\hat{p}_5^2}{\hat{p}_2}} \text{,}
\end{eqnarray*}
 where $\hat{p}_2$ is given by the quartic equation
  \begin{dmath*}
  -40\hat{p}_2^4 +(76w-92u+16v)\hat{p}_2^3 + \big[ -4u^2 +4(u-w)(3w-5u+2v)-2(w+v-u)^2\big] \hat{p}_2^2 + \big[(w-u)(w+v-u)-2u^2(3u+v-w) \big] \hat{p}_2 + u^4 = 0 \text{ , } \end{dmath*}
where $\hat{p}_3$ is given by the quartic equation
 \begin{dmath*}
  10 \hat{p}_3^4 -2(10u+3v+2w)\hat{p}_3^3 + \big[ 6u^2 +(4u+v)(3u+v+w)\big] \hat{p}_3^2 - u^2(7u+2v+w) \hat{p}_3 + u^4 = 0 \text{ , } \end{dmath*}
and $\hat{p}_5$ is given by the quartic equation
  \begin{dmath*}
 5 \hat{p}_5^4 +(6v+4w)\hat{p}_5^3 + \big[ 6u^2 +(2v-4u)(v+w)\big] \hat{p}_5^2 + \big[-4u^3-2v (u^2+2uv+2wu) \big] \hat{p}_5 + \big[ u^4+2uv(uv+wu) \big] = 0 \end{dmath*}
  with $u=\frac{2u_3+u_5}{N}$, $v=\frac{u_1-u_4}{N}$ and $w=\frac{2u_2+u_4+u_5}{N}$.

\subsection{Toric Del Pezzo surfaces of degree greater than five.}
When the degree of the Del Pezzo surface is greater than or equal to five, each of the defining equations of the probability distribution $(\hat{p}_1, \hat{p}_2, \hat{p}_3, \cdots \hat{p}_n) $ satisfy an equation of degree five or higher. By the Abel-Ruffini theorem there is no algebraic solution for a general polynomial equation of degree five or higher, therefore one cannot obtain a closed form solution for the maximum likelihood estimate in these cases.

Also, Table \ref{table:1} provides a list of the ML degrees of toric Del Pezzo surfaces of degree up to six using Macaulay2 and Algorithm 6 in Section 4 of \cite{Hosten_etal}.
\begin{table}[h]
\centering
\begin{tabular}{||c c||}
 \hline
 Ideal of some degree $d$ Del Pezzo $S_d$ & ML-degree\\ [0.5ex]
 \hline\hline
 $S_3: p_1 p_2 p_3 - p_4^3$ & 3  \\
 \hline
 $S_{4}: p_1p_4-p_5^2, p_2p_3-p_1p_4$ & 4 \\
 \hline
 $S_{4}^{'}: p_2p_4-p_3p_5, p_1p_3-p_5^2$ & 4 \\
 \hline
 $S_{4}^{''}: p_2p_4-p_3^2, p_1p_3-p_5^2$ & 4 \\
 \hline
$S_{5}: p_3p_5-p_4p_6, p_2p_5-p_6^2, p_2p_4-p_3p_6$ & \\
$p_1p_4-p_6^2, p_1p_3-p_2p_6$ & 3\\
\hline
$S_{5}^{'}: p_3p_5-p_4p_6, p_2p_5-p_6^2, p_2p_4-p_3p_6$ & \\
$ p_1p_4-p_2p_6, p_2^2-p_1p_3$ & 5 \\
\hline
$S_{6}: p_4p_6-p_5p_7, p_3p_6-p_7^2, p_2p_6-p_1p_7$ & \\
$p_3p_5-p_4p_7, p_2p_5-p_7^2, p_1p_5-p_6p_7$ & \\
$p_2p_4-p_3p_7, p_1p_4-p_7^2, p_1p_3-p_2p_7$ & 6 \\
\hline
$S_{6}^{'}: p_5 p_6-p_1 p_7, p_4 p_6-p_2 p_7, p_3 p_5-p_4 p_7$ & \\
$ p_2 p_5-p_7^2, p_2 p_4-p_3 p_7, p_1 p_4-p_7^2$ & \\
$ p_1 p_3-p_2 p_7, p_2^2-p_3 p_6, p_1 p_2-p_6 p_7$ & 6 \\
\hline
$S_{6}^{''}: p_6^2-p_5 p_7, p_4 p_6-p_3 p_7, p_3 p_6-p_2 p_7$ & \\
$ p_4 p_5-p_2 p_7, p_3 p_5-p_2 p_6, p_2 p_4-p_1 p_7$ & \\
$  p_3^2-p_1 p_7,p_2 p_3-p_1 p_6, p_2^2-p_1 p_5$ & 6 \\
\hline
$S_{6}^{'''}: p_6^2-p_5p_7, p_5p_6-p_4p_7, p_3p_6-p_2p_7$ & \\
$ p_5^2-p_4p_6, p_3p_5-p_2p_6, p_3p_4-p_2p_5$ & \\
$ p_3^2-p_1 p_6, p_2 p_3-p_1 p_5, p_2^2-p_1p_4$ & 6 \\ [1ex]
\hline
\end{tabular}
\caption{ML degree of some Del Pezzo surfaces}
\label{table:1}
\end{table}

 Already for the case of a septic Del Pezzo surface corresponding to the reflexive polytope $7a$ in Table~\ref{table:2}, the memory of a Dell laptop with 2.4 GHz Intel Core i5 processor was not sufficient for the computation to be completed.  However, the results described in Table~\ref{table:1} are in accordance with \cite[Theorem 3.2]{Huh_Sturmfels}, which states that the ML degree of a projective toric variety is bounded above by its degree. The cubic and quartic Del Pezzo surfaces in Table~\ref{table:1} are complete intersections and have ML degree which equals the degree of the Del Pezzo surface. It is known, that Del Pezzo surfaces of degree greater than or equal to five are no longer complete intersections (see \cite{Dolgachev}). We see in Table~\ref{table:1} that the ML degree drops to three in the case of a quintic Del Pezzo surface $S_5$ corresponding to the reflexive polytope $5b$ in  Table~\ref{table:2}.

The natural next step is to consider the maximum likelihood estimation problem for scaled toric Del Pezzo surfaces which give different embeddings of isomorphic Del Pezzo surfaces in projective space and compute their ML degrees, as is the line of research in \cite{Carlos_etal}.

The author would like to thank the referees for their useful comments and suggestions which greatly improved this paper.

\noindent
Dimitra Kosta\\
School of Mathematics and Statistics,\\
University of Glasgow,\\
Glasgow G12 8SQ, UK.\\
Dimitra.Kosta@glasgow.ac.uk


\addcontentsline{toc}{section}{Bibliography}
        \begin{thebibliography}{99}
        \label{Bibliography}

\bibitem{Carlos_etal}
C.\, Amendola, N.\, Bliss, I.\, Burke, C.\,R.\, Gibbons, M.\, Helmer, S.\, Hosten, E.\, D.\, Nash, J.\, I.\, Rodriguez, D.\, Smolkin, \emph{The maximum likelihood degree of toric varieties}, https://arxiv.org/pdf/1703.02251.pdf.

\bibitem{Bishop_etal}
Y.\, M.\, Bishop, S.\,E.\, Fienberg and P.\,W.\, Holland, \emph{Discrete Multivariate Analysis: Theory and Practice}, Springer-Verlag (2007), New York.

         \bibitem{Catanese_etal}
         F.\,Catanese, S.\,Hos̡̡ten, A.\,Khetan and  B.\,Sturmfels, \emph{The Maximum likelihood degree}, American Journal of Mathematics, \textbf{18} (2006), 671-697.

\bibitem{Dolgachev}
I.\, V.\, Dolgachev, \emph{Del Pezzo surfaces} in Classical Algebraic Geometry: A Modern View, Cambridge University Press  (2012), 347-425.

\bibitem{Drton_etal}
M.\,Drton, B.\,Sturmfels, S.\,Sullivant, \emph{Lectures on Algebraic Statistics}, Oberwolfach Seminars, \textbf{39} (2008), Birkhäuser, Basel.

\bibitem{Hauenstein_etal}
J.\,D.\,Hauenstein, J.\,Rodriguez, and B.\,Sturmfels, \emph{Maximum likelihood for matrices with rank constraints}, Journal of Algebraic Statistics, \textbf{5} (2014), 18-38.

\bibitem{Hosten_etal}
S.\,Hos̡̡ten, A.\,Khetan and B.\,Sturmfels, \emph{Solving the likelihood equations}, Foundations of Computational Mathematics, \textbf{5} (2005), 389-407.

\bibitem{Huh}
J.\, Huh, \emph{Varieties with maximum likelihood degree one}, Journal of Algebraic Statistics, \textbf{5} (2014), 1-17.

\bibitem{Huh_Sturmfels}
J.\, Huh, B.\, Sturmfels, Likelihood Geometry in \emph{Combinatorial Algebraic Geometry}, Lecture Notes in Mathematics, \text{2108} (2014), 63-117.


\bibitem{Kollar}
J.\,Koll\'ar, \emph{Minimal models of algebraic threefolds: Mori's program}, Astérisque \text{177} (1989), 303-326.


\bibitem{Kreuzer_Skarke}
M.\,Kreuzer, H.\,Skarke, \emph{ On the classification of reflexive polyhedra}, Commun. Math. Phys. \text{185} (1997) 495-508.

\bibitem{Lin_etal}
S.\, Lin, B.\, Sturmfels, and Z.\, Xu, \emph{Marginal Likelihood Integrals for Mixtures of Independence Models}, Journal of Machine Learning Research, \textbf{10} (2009), 1611-1631.


\bibitem{Pachter_Sturmfels}
L.\,Pachter and B.\,Sturmfels, \emph{Algebraic Statistics for Computational Biology}, (2005) Cambridge University Press, New York.

\bibitem{Reid}
M.\, Reid, \emph{Canonical 3-folds}, Journ\'ees de G\'eometrie Alg\'ebrique d’Angers, Juillet 1979/Algebraic Geometry, Angers, 1979, Sijthoff \& Noordhoff, Alphen aan den Rijn (1980), 273-310.

\bibitem{Sato}
H.\,Sato, \emph{Toward the classification of higher-dimensional Fano varieties}, Tohoku Math. J., \textbf{52} (2000), 383-413.

\bibitem{Sturmfels2009}
B.\,Sturmfels, Open problems in Algebraic Statistics, In \emph{"Emerging Applications of Algebraic Geometry"}, I.M.A. Volumes in Mathematics and its Applications, \textbf{149} (2009), Springer, New York, 351-364.

\bibitem{Wisniewski_Buczynska_2007}
W. Buczynska, J. Wisniewski, \emph{On the geometry of binary symmetric models of phylogenetic trees}, J. Euro. Math. Society, \textbf{9} (2007), no. 3, 609-635.

  	\end{thebibliography}
\end{document}